\documentclass[12pt]{amsart}
\usepackage{latexsym}
\usepackage{amsmath}
\usepackage{amsfonts}
\usepackage{amssymb}
\usepackage{amsthm}

\theoremstyle{plain}
\newtheorem{theorem}{Theorem}[section]
\newtheorem{lemma}[theorem]{Lemma}
\newtheorem{proposition}[theorem]{Proposition}

\newtheorem*{theorem*}{Theorem}

\theoremstyle{definition}

\newtheorem{example}{Example}

 \theoremstyle{remark}

\font\tenmsb=msbm10 at 12pt \font\sevenmsb=msbm7 at 8pt
\font\fivemsb=msbm5 at 6pt
\newfam\msbfam
\textfont\msbfam=\tenmsb \scriptfont\msbfam=\sevenmsb
\scriptscriptfont\msbfam=\fivemsb

\newcommand{\bbE}{\mathbb{E}}
\newcommand{\bbN}{\mathbb{N}}

\newcommand{\bbR}{\mathbb{R}}
\newcommand{\bbT}{\mathbb{T}}
\newcommand{\bbZ}{\mathbb{Z}}

\newcommand{\cG}{\mathcal{G}}

\newcommand{\cX}{\mathcal{X}}
\newcommand{\cY}{\mathcal{Y}}

\newcommand{\cZ}{\mathcal{Z}}

\newcommand{\gG}{\Gamma}

\begin{document}
\title{Polynomial averages converge to the product of integrals}

\author{Nikos Frantzikinakis and Bryna Kra}

\address{Department of Mathematics, McAllister Building,
The Pennsylvania State University,
University Park, PA 16802}
\email{nikos@math.psu.edu}
\email{kra@math.psu.edu}

\begin{abstract}
We answer a question posed by Vitaly Bergelson, showing that in a
totally ergodic system, the average of a product of functions
evaluated along polynomial times, with polynomials of pairwise
differing degrees, converges in $L^{2}$ to the product of the
integrals.  Such averages are characterized by nilsystems and so
we reduce the problem to one of uniform distribution of polynomial
sequences on nilmanifolds.
\end{abstract}
\maketitle

\section{Introduction}

\subsection{Bergelson's Question}

In \cite{Bergelson}, Bergelson asked if the average of a product
of functions in a totally ergodic system (meaning that
each power of the
transformation is ergodic) evaluated along
polynomial times converges in $L^{2}$ to the product of the
integrals. More precisely, if $(X, \cX, \mu, T)$ is a totally
ergodic probability measure preserving system, $p_1, p_{2},
\ldots, p _{k}$ are polynomials taking integer values on the
integers with pairwise distinct non-zero degrees, and $f_1,
f_{2},\ldots, f_k \in L^\infty(\mu)$, does
$$
\lim_{N\to\infty}\Big\Vert \frac{1}{N}\sum_{n=0}^{N-1}
f_1(T^{p_1(n)}x)f_2(T^{p_2(n)}x)\ldots f_k(T^{p_k(n)}x) -
\prod_{i=1}^k\int f_i\,d\mu \Big\Vert_{L^2(\mu)}
$$
equal $0$?

We show that the answer to this question is positive, under
slightly more general assumptions.  We start with some definitions
in order to precisely state the theorem.

An  \emph{integer polynomial} is a polynomial taking integer
values on the integers. A family of integer  polynomials $
\{p_1(n), p_{2}(n), \ldots, p_k(n)\} $
 is said to be \emph{independent} if for all integers $m_1,
 m_{2}, \ldots, m_{k}$ with at least some $m_{j}\neq 0 $, $j
 \in\{1, 2, \ldots, k\}$,
the polynomial $\sum_{j=1}^k m_jp_j(n)$ is not constant.

We prove:
\begin{theorem}
\label{T:product}
Let $(X, \cX, \mu, T)$ be a totally ergodic
measure preserving probability system and assume that $\{p_1(n),
p_{2}(n), \ldots, p_k(n)\}$ is an independent family of
polynomials. Then for $f_1, f_{2}, \ldots, f_k \in L^\infty(\mu)$,
\begin{equation}
\label{eq:polynomial} \lim_{N\to\infty}\Big\Vert
\frac{1}{N}\sum_{n=0}^{N-1} f_1(T^{p_1(n)}x)f_2(T^{p_2(n)}x)\ldots
f_k(T^{p_k(n)}x) - \prod_{i=1}^k\int f_i\,d\mu
\Big\Vert_{L^2(\mu)}
\end{equation}
equals $0$.
\end{theorem}

The assumption that the polynomial family is independent is
necessary, as can be seen by considering an irrational rotation on
the circle.  An ergodic rotation on a finite group with at least
two  elements demonstrates that the hypothesis of total ergodicity
is necessary; in this example, the average for any independent
family with $k>1$ polynomials does not converge to the product of
the integrals for appropriate choice of the functions $f_i$.

If one assumes that $T$ is weakly mixing, Bergelson
\cite{Berg3} showed that for all polynomial families, the limit
in~\eqref{eq:polynomial} exists and is constant. However, without
the assumption of weak mixing one can easily show that the limit
need not be constant, even when restricting to polynomials of
degree one. For the polynomial families $(n, n^2)$ and $(n^2,
n^2+n)$, the convergence to the product of the integrals was
proved by Furstenberg and Weiss \cite{FW}.   The existence of the
limit in a totally ergodic system
for an arbitrary family of integer polynomials was shown in
Host and Kra \cite{HK}, but further analysis is needed to describe
the form of the limit.

\subsection{Reduction to a problem of uniform distribution}

In \cite{HK}, Host and Kra showed that for any family of
polynomials, the characteristic factor of the average
in~\eqref{eq:polynomial} in a totally ergodic system
is an inverse limit of nilsystems. We
need a few definitions to make this statement precise.

Given a group $G$, we denote the commutator of $g,h\in G$
by $[g,h]=g^{-1}h^{-1}gh$.  If $A, B\subset G$, then
$[A,B]$ is defined to be $\{[a,b]:a\in A, b\in B\}$.
A group
$G$ is said to be {\em $k$-step nilpotent} if its $(k+1)$
commutator $[G, G^{(k)}]$ is trivial.  If $G$ is a $k$-step
nilpotent Lie group and $\gG$ is a discrete cocompact subgroup,
then the compact space $X = G/\gG$ is said to be a {\em $k$-step
nilmanifold}.  The group $G$ acts on $G/\gG$ by left translation and
the translation by a fixed element
$a\in G$ is given by $T_{a}(g\gG) = (ag) \gG$.  Let
$\mu$ denote the unique probability measure on $X$ that is
invariant under the action of $G$ by left translations (called the
{\em Haar measure}) and let $\cG/\gG$ denote the Borel $\sigma$-algebra of
$G/\gG$.  Fixing an element $a\in G$, we call
the system $(G/\gG, \cG/\gG, \mu, T_{a})$ a {\em $k$-step
nilsystem} and call the map $T_a$ a {\em nilrotation}.

A {\em factor} of the measure preserving system $(X,\cX, \mu,T)$
is a measure preserving system $(Y,\cY,\nu,S)$ so that there
exists a measure preserving map $\pi:X\to Y$ taking $\mu$ to $\nu$
and such that $S\circ \pi = \pi \circ T$.  In a slight abuse of
terminology, when the underlying measure space is implicit we call
$S$ a factor of $T$.

In this terminology, Host and Kra's result means that there exists
a factor $(Z,\cZ,m)$ of $X$, where $\cZ$ denotes the Borel
$\sigma$-algebra of $Z$ and $m$ its Haar measure, so that the action of
$T$ on $Z$ is an inverse limit of nilsystems and furthermore, whenever
$\bbE(f_j|\cZ) = 0 $ for some $j \in \{1, 2, \ldots, k\}$, the
average in~\eqref{eq:polynomial} is itself $0$. Since an inverse
limits of nilsystems can be approximated arbitrarily well by a
nilsystem, it suffices to verify Theorem~\ref{T:product} for
nilsystems.  Moreover, since measurable functions can be
approximated arbitrarily well in $L^2$ by  continuous functions,
Theorem~\ref{T:product} is equivalent to  the following
generalization of Weyl's polynomial uniform distribution theorem (see
Section~\ref{sec:uniform} for the statement of Weyl's Theorem):

\begin{theorem}\label{T:equi}
Let $X=G/\gG$ be a nilmanifold, $(G/\gG, \cG/\gG, \mu, T_{a})$ a
nilsystem and suppose that the nilrotation $T_a$ is totally
ergodic. If $\{p_1(n), p_2(n), \ldots, p_k(n)\}$ is an independent polynomial
family, then for almost every $x\in X$ the sequence $(a^{p_1(n)}x,
a^{p_{2}(n)}x, \ldots,a^{p_k(n)}x)$ is uniformly distributed in
$X^k$.
\end{theorem}

If $G$ is connected, we can reduce Theorem~\ref{T:equi} to a
uniform distribution problem that is easily verified using the
standard uniform distribution theorem of Weyl. The general (not
necessarily connected) case is more subtle. Using a result of
Leibman \cite{Leibman}, in Section~\ref{sec:prelim},
we reduce the problem to studying the action
of a polynomial sequence on a factor space with abelian identity
component. The key step (Section~\ref{sec:affine})
is then to prove that nilrotations acting
on such spaces are isomorphic to affine transformations on some
finite dimensional torus.  In Section~\ref{sec:uniform}, we
complete the proof by checking
the result for affine transformations.

\section{Reduction to an abelian connected component}
\label{sec:prelim}

Suppose that $G$ is a nilpotent Lie group and $\Gamma$ is a
discrete, cocompact subgroup.  Throughout, we let $G_0$ denote the
connected component of the identity element and denote the
identity element by $e$.

A sequence $g(n) = a_1^{p_1(n)}a_2^{p_2(n)}\ldots
a_k^{p_k(n)}$ with $a_1, a_{2}, \ldots, a_k \in G$ and $p_1,
p_{2}, \ldots, p_k$ integer polynomials is called a
\emph{polynomial sequence} in $G$. We are interested in studying
uniform distribution properties of polynomial sequences on the
nilmanifold $X=G/\gG$.

Leibman \cite{Leibman} showed that the uniform distribution of a
polynomial sequence in a connected nilmanifold reduces to uniform distribution
in a certain factor:

\begin{theorem*}{\bf [Leibman]}
Let $X = G/\gG$ be a connected nilmanifold and let $g(n) =
a_1^{p_1(n)}a_2^{p_2(n)}\ldots a_k^{p_k(n)}$ be a polynomial
sequence in $G$. Let $Z=X/[G_0,G_0]$  and let $\pi\colon X\to Z$
be the natural projection. If $x\in X$ then $\{g(n)x\}_{n\in\bbZ}$
is uniformly distributed in $X$ if and only if
$\{g(n)\pi(x)\}_{n\in\bbZ}$ is uniformly distributed in $Z$.
\end{theorem*}

We remark that if $G$ is connected, then the factor $X/[G_0,G_0]$ is an
abelian group.  However, this does not hold in general as the
following examples illustrate:

\begin{example}
\label{ex:one}
On the space $G=\bbZ\times\bbR^2$,
define multiplication as follows: \\
if $g_1=(m_1,x_1,x_2)$ and $g_2=(n_1,y_1,y_2)$, let
$$
g_1\cdot g_2=(m_1+n_1,x_1+y_1, x_2+y_2+m_1y_1).
$$ Then $G$ is a $2$-step
nilpotent group and $G_0=\{0\}\times \bbR^2$ is abelian. The
discrete subgroup $\gG=\bbZ^3$ is  cocompact and $X=G/\gG$ is
connected. Moreover,  $[G_0,G_0]=\{{\bf e}\}$ and so
$X/[G_0,G_0]=X$.
\end{example}

\begin{example}
On the space $G=\bbZ\times \bbR^3$,
define multiplication as follows: \\
 if $g_1=(m_1,x_1,x_2,x_3)$ and
$g_2=(n_1,y_1,y_2,y_3)$, let
$$
g_1\cdot g_2=( m_1+n_1, x_1+y_1, x_2+y_2+m_1y_1,
x_3+y_3+m_1y_2+\frac{1}{2}m_1^2y_1).
$$
Then $G$ is a $3$-step nilpotent group and $G_0=\{0\}\times\bbR^3$
is abelian. The discrete subgroup $\gG=\bbZ^3\times (\bbZ/2)$ is
cocompact and $X=G/\gG$ is connected. Again, $X/[G_0,G_0]=X$.
\end{example}

We use Leibman's theorem to reduce the problem on uniform
distribution to the case that $G_0$ is abelian:
\begin{proposition}\label{P:red1}  Theorem~\ref{T:equi} follows if
it holds for all nilsystems $(G/\gG,\cG/\gG,\mu,T_a)$ with $G_0$
abelian and $T_a$ totally ergodic.
\end{proposition}

\begin{proof}
 Given $a\in G$ and $x\in X=G/\gG$, let $a_1=(a,e,\ldots,e),
a_2=(e, a, e, \ldots, e), \ldots, a_k=(e,e,\ldots,a)\in G^k$,
$\tilde{x}=(x,\ldots,x)\in X^k$, and
$g(n)=T_{a_1}^{p_1(n)}T_{a_2}^{p_2(n)} \cdots T_{a_k}^{p_k(n)}$.
We need to check that for $\mu$-a.e. $x\in X$ the polynomial
sequence $g(n)\tilde{x}$ is uniformly distributed in $X^k$. By
Leibman's Theorem,  it suffices to check that $g(n)\pi(\tilde{x})$
is uniformly distributed in the nilmanifold $Z^k$, where
$Z=X/[G_0,G_0]$ and $\pi\colon G\to G/[G_0,G_0]$ is the natural
projection. Since $(G/[G_0,G_0])_0$ is abelian and a factor of a
totally ergodic system is totally ergodic, the statement follows.
\end{proof}

\section{Reduction to an affine transformation on a torus}
\label{sec:affine} We reduce the problem on uniform distribution
(Theorem~\ref{T:equi}) to studying an affine transformation on a
torus. If $G$ is a group then a map $T\colon G\to G$ is said to be
\emph{affine} if $T(g) = bA(g)$ for a homomorphism $A$ of $G$ and
some $b\in G$.  The homomorphism $A$ is said to be {\em unipotent}
if there exists $n\in\bbN$ so that so that $(A-{\text Id})^{n}=0$.
In this case we say that the affine transformation $T$ is a
unipotent affine transformation.

\begin{proposition}\label{C:1}
Let $X=G/\gG$ be a connected nilmanifold such that $G_0$ is
abelian. Then any nilrotation $T_a(x)=ax$ defined on $X$ with the
Haar measure $\mu$ is isomorphic to a unipotent affine
transformation on some finite dimensional torus.
\end{proposition}
\begin{proof}
First observe that for every $g\in G$, the subgroup $g^{-1}G_0 g$
is both open and closed in $G$ so $g^{-1}G_0 g=G_0$. Hence, $G_0$ is a
normal subgroup of $G$. Similarly, since $G_0\gG$ is both open and closed
in $G$, we have that $(G_0\gG)/\gG$ is open and closed in $X$.
Since $X$ is connected, $X=(G_0\gG)/\gG$ and so
$G=G_0\gG$.

We claim that $\gG_0=\gG\cap G_0$ is a normal subgroup of $G$. Let
$\gamma_0\in \Gamma_0$ and  $g=g_0\gamma$, where $g_0\in G_0$ and
$\gamma\in \Gamma$. Since $G_0$ is normal in $G$, we have that
$g^{-1}\gamma_0 g\in G_0$. Moreover,
$$
g^{-1} \gamma_0 g=\gamma^{-1} g_0^{-1} \gamma_0 g_0 \gamma=
\gamma^{-1} \gamma_0 \gamma \in\gG,
$$
the last equality being valid since $G_0$ is abelian. Hence,
$g^{-1}\gamma_0 g\in\Gamma_0$ and $\Gamma_0$ is normal in $G$.

Therefore we can substitute $G/\gG_0$ for $G$ and $\gG/\gG_0$ for
$\gG$; then $X=(G/\gG_0)/(\gG/\gG_0)$. So we can assume that
$G_0\cap \gG=\{e\}$. Note that we now have that $G_0$ is a
connected compact abelian Lie group and so is isomorphic to some
finite dimensional torus $\bbT^d$.

Every $g\in G$ is uniquely representable in the form
$g=g_0\gamma$, with $g_0\in G_0$, $\gamma \in \gG$. The map
$\phi\colon X\to G_0$, given by $\phi(g\Gamma)=g_0$ is a well
defined homeomorphism. Since   $\phi(hg\Gamma)=h\phi(g\Gamma)$ for
any $h\in G_0$, the measure $\phi(\mu)$ on $G_0$ is invariant under left
translations.  Thus $\phi(\mu)$ is the Haar measure on $G_0$. If
$a=a_0\gamma$, $g=g_0\gamma'$ with $a_0,g_0\in G_0$ and
$\gamma,\gamma'\in\gG$, then $ag\gG=a_0\gamma g_0 \gamma^{-1}\gG$.
Since $\gamma g_0 \gamma^{-1} \in G_0$, we have that
$\phi(ag\Gamma)=a_0\gamma g_0\gamma^{-1}$. Hence $\phi$ conjugates
$T_a$ to $T_a'\colon G_0\to G_0$ defined by
$$
T_a'(g_0)=\phi T_a \phi^{-1}=a_0\gamma g_0 \gamma^{-1}.
$$
 Since $G_0$ is abelian
this is an affine map; its linear part $g_0\mapsto \gamma g_0
\gamma^{-1}$ is unipotent since $G$ is nilpotent.  Letting
$\psi\colon G_0\to \bbT^d$ denote the isomorphism between $G_0$ and
$\bbT^d$, we have that $T_a$ is isomorphic to the unipotent affine
transformation $S=\psi T_a'\psi^{-1}$ acting on $\bbT^d$.
\end{proof}

We illustrate this with the examples of the previous section:
\begin{example}
Let $X$ be as in Example $1$ and let $a=(m_1,a_1,a_2)$. Since
$G_0/\gG_0=\bbT^2$ we see that  $T_a$ is isomorphic to the
unipotent affine
transformation $S \colon \bbT^2\to\bbT^2$ given by
$$
S(x_1,x_2)=(x_1+a_1,x_2+m_1x_1+a_2).
$$
\end{example}

\begin{example} Let $X$ be as in Example $2$ and
$a=(m_1,a_1,a_2,a_3)$. Since $G_0/\gG_0=\bbR^3/(\bbZ^2
\times\bbZ/2)$, and $\psi\colon G_0/\gG_0\to \bbT^3$ defined by
$\psi(x_1,x_2,x_3)=(x_1,x_2,2x_3)$ is an isomorphism, we see that
$T_a$ is isomorphic to the unipotent affine transformation $S \colon
\bbT^3\to\bbT^3$ given by
$$
S(x_1,x_2,x_3)=(x_1+a_1,x_2+m_1 x_1+ a_2,x_3+
2m_1x_2+m_1^2x_1+2a_3).
$$
\end{example}

\begin{proposition}\label{P:red2}  Theorem~\ref{T:equi} follows if
it holds for all  nilsystems $(G/\gG,\cG/\gG,\mu,T_a)$ such that
$T_a$ is isomorphic to an ergodic, unipotent, affine
transformation on some finite dimensional torus.
\end{proposition}
\begin{proof}
We first note that since $X=G/\gG$ admits a totally ergodic
nilrotation $T_a$, it must be connected. Indeed, let $X_0$ be the
identity component of $X$. Since $X$ is compact, it is a disjoint
union of $d$ copies of translations of $X_0$ for some $d\in\bbN$.
Since $a$ permutes these copies, $a^d$ preserves $X_0$. By
assumption the translation by  $T_{a^d}=T_a^d$ is ergodic and so
$X_0=X$.

By Proposition~\ref{P:red1} we can assume that $G_0$ is abelian.
Since $X$ is connected, the result follows from
Proposition~\ref{C:1}.
\end{proof}

\section{Uniform distribution for an affine transformation}
\label{sec:uniform}

We are left with showing that  Theorem~\ref{T:equi} holds when the
nilsystem is isomorphic to an ergodic,  unipotent,  affine system
on a finite dimensional torus. Before turning into the proof,
note that if $G$ is connected then the uniform distribution
property of Theorem~\ref{T:equi} holds for every $x\in X$.
However, this does not hold in general. We illustrate this with
the following example:

\smallskip

\begin{example}We have seen that the nilrotation of
Example 1 is isomorphic to the affine transformation $S \colon
\bbT^2\to\bbT^2$ given by
$$
S(x_1,x_2)=(x_1+a_1,x_2+m_1x_1+a_2).
$$
If $m_1=2$ and $a_1=a_2=a$  is irrational then  $S$ is totally
ergodic and $S^n(x_1,x_2)=(x_1+na,x_2+2nx_1+n^2a)$. Then
$$
\bigl(S^n(0,0),S^{n^2}(0,0)\bigr)=(na, n^2a,n^2a,n^4a)
$$
is not uniformly distributed on $\bbT^4$. On the other hand
\begin{eqnarray*}
\lefteqn{\bigl(S^n(x_1,x_2),S^{n^2}(x_1,x_2)\bigr)= }\\
& & (x_1+na, x_2+2nx_1+n^2a,x_1+n^2a,  x_2+2n^2x_1+n^4a,)
\end{eqnarray*}
is uniformly distributed on $\bbT^4$ as long as $a$ and $x_1$ are
rationally independent.
\end{example}

The main tool used in the proof of Theorem~\ref{T:equi} is  the following
classic theorem of Weyl \cite{W} on uniform distribution:
\begin{theorem*}{\bf [Weyl]}\label{T:Weyl}
(i) Let $a_n\in \bbR^d$. Then $a_n$ is uniformly distributed in
$\bbT^d$ if and only if
$$
\lim_{N\to\infty}\frac{1}{N}\sum_{n=1}^N e^{2\pi i m\cdot a_n}=0
$$
for every nonzero $m\in\bbZ^d$, where $m\cdot a_n$ denotes the
inner product of $m$ and $a_n$.

(ii) If \ $a_n=p(n)$ where $p$ is a real valued polynomial with at
least one nonconstant coefficient irrational then
$$
\lim_{N\to\infty}\frac{1}{N}\sum_{n=1}^N e^{2\pi i a_n}=0.
$$
\end{theorem*}

Before turning to the proof of Theorem~\ref{T:equi}, we prove
a lemma that simplifies the computations:
\begin{lemma} \label{L:simplify}
Let $T\colon\bbT^d\to \bbT^d$ be defined by $T(x)=Ax+b$, where $A$
is a $d\times d$ unipotent integer matrix and $b\in\bbT^d$. Assume
furthermore that $T$ is ergodic. Then $T$ is a factor of an
ergodic affine transformation $S\colon\bbT^d\to \bbT^d$, where
$S=S_1\times S_{2}\times\cdots\times S_s$ and for $r=1, 2, \ldots,
s$,
$S_r\colon\bbT^{d_r}\to \bbT^{d_r}$  ($\sum_{r=1}^sd_r=d$) has the
form
$$
S_r(x_{r1}, x_{r2},\ldots,x_{rd_r})=(x_{r1}+b_{r},x_{r2}+
x_{r1},\ldots,x_{rd_r}+x_{rd_r-1})
$$
for some $b_r\in\bbT$.
\end{lemma}
\begin{proof}
Let $J$ be the Jordan canonical form of $A$ with Jordan blocks
$J_r$ of dimension $d_r$ for $r=1,2,\ldots,s$. Since $A$ is
unipotent, all diagonal entries of $J$ are equal to $1$. There
exists a matrix $P$ with rational entries such that $PA=JP$. After
multiplying $P$ by an appropriate integer, we can assume that it
too has integer entries. So  $P$ defines a homomorphism
$P\colon\bbT^d\to \bbT^d$ such that $PT=SP$, where
$S\colon\bbT^d\to \bbT^d$ is given by $S(x)=J(x)+c$ for $c=P(b)$.
Hence, $T$ is a factor of $S$.  By making the change of variables
$x_{ij}\to x_{ij}+a_{ij}$, we can assume that $S$ has the
advertised form.

It remains to show that $S$ is ergodic. Since $J$ is
unipotent, using a theorem of Hahn (\cite{H}, Theorem 4) we get that
ergodicity of $S$ is equivalent  to showing that for every
nontrivial character $\chi$ in the dual of $\bbT^d$ we have the
implication
$$
\chi(Jx)=\chi(x) \text{ for every } x\in\bbT^d \Rightarrow
\chi(c)\neq 1.
$$
Suppose that $ \chi(Jx)=\chi(x)$. Using the relation $PA=JP$ we
get that $\chi'(Ax)=\chi'(x)$ where $\chi'(x)=\chi(Px)$. Since
$T(x)=Ax+b$ is assumed to be ergodic, again using Hahn's theorem
we get that $ \chi'(b)\neq 1$. The relation $PA=JP$ implies that $
\chi(c)\neq 1 $ and the proof is complete.
\end{proof}
\begin{proof}[Proof of Theorem~\ref{T:equi}] By Proposition \ref{P:red2}
it suffices to verify the uniform distribution property for all
ergodic, unipotent, affine transformations on $\bbT^d$. First
observe that relation \eqref{eq:polynomial} of
Theorem~\ref{T:product} is preserved when passing to factors.
Hence, using Lemma~\ref{L:simplify} we can assume that
$T=T_1\times T_{2}\times\cdots\times T_s$, where
$T_r\colon\bbT^{d_r}\to \bbT^{d_r}$ ($\sum_{r=1}^sd_r=d$) is given
by
$$
T_r(x_{r1}, x_{r2},\ldots,x_{rd_r})=(x_{r1}+b_{r},x_{r2}+
x_{r1},\ldots,x_{rd_r}+x_{rd_r-1}),
$$
for $r=1, 2, \ldots,s$.
 Since $T$ is ergodic the set $\{b_1, b_2, \ldots, b_s\}$ is
rationally independent. For convenience, set $x_{r0}=b_r$ for
$r=1, 2,\ldots s$.

We claim that if $x$ is chosen so that the set $A=\{x_{rj}: 1\leq
r \leq s, 0\leq j \leq d_r\}$ is rationally independent, then the
polynomial sequence $g(n)\tilde{x}=$ $(T^{p_1(n)}x,
T^{p_2(n)}x,\ldots, T^{p_k(n)}x)$ is uniformly distributed on
$\bbT^{dk}$ (we include $x_{rd_r}$ in $A$ only for simplicity). To
see this we use the first part of Weyl's theorem; letting
$Q_{rjl}(n)$ denote the $j$-th coordinate of $T_r^{p_l(n)}x$ and
\begin{equation}\label{E:R(n)}
R(n)=\sum_{r,j,l} m_{rjl} Q_{rjl}(n)
\end{equation}
where $\{m_{rjl}: 1 \leq r \leq s, 1\leq j \leq d_r, 1\leq l \leq
k\}$ are integers, not all of them zero, it suffices to check that
\begin{equation}\label{E:Weyl}
\lim_{N\to\infty}\frac{1}{N}\sum_{n=1}^N e^{2\pi i R(n)}=0.
\end{equation}
To prove \eqref{E:Weyl} we use the second part of Weyl's theorem;
it suffices to show that the polynomial $R(n)$ has at least one
nonconstant coefficient irrational.
 We compute
\begin{equation}\label{E:Q}
Q_{rjl}(n)=x_{rj}+\binom{p_l(n)}{1}x_{rj-1}+\cdots+
\binom{p_l(n)}{j-1} x_{r1}+\binom{p_l(n)}{j}x_{r0}.
\end{equation}
We can put $R(n)$ in the form
\begin{equation}\label{E:R}
R(n)=\sum_{r,j}R_{rj}(n)x_{rj},
\end{equation}
where $R_{rj}$ are integer polynomials and $1\leq r\leq s$, $0\leq
j \leq d_r$. This representation is unique since the $x_{rj}$ are
rationally independent. So it remains to show that some $R_{rj}$
is nonconstant. To see this, choose any $r_0$ such that
$m_{r_0jl}\neq 0$ for some $j,l$, and define $j_0$ to be the
maximum $1\leq j \leq d_{r_0}$ such that $m_{r_0jl}\neq 0$ for
some $1\leq l \leq k$. We show that $R_{r_0,j_0-1}$ is
nonconstant. By the definition of  $j_0$ we have $m_{r_0jl}=0$ for
$j>j_0$. For $j\leq j_0$ we see from \eqref{E:Q} that  the
variable $x_{r_0j_0-1}$ appears only in the polynomials
$Q_{r_0j_0l}$ with coefficient $p_l(n)$, and if $j_0>1$ also in
the polynomials $Q_{r_0(j_0-1)l}$ with coefficient $1$. It follows
from \eqref{E:R(n)}  and \eqref{E:R} that
$$
R_{r_0j_0-1}(n)=\sum_{l=1}^k m_{r_0j_0l} p_l(n) +c,
$$
where $c=\sum_{l=1}^k m_{r_0j_0l}$ if $j_0>1$, and $c=0$ if
$j_0=1$. Since the polynomial family $\{p_i(n)\}_{i=1}^k$ is
independent and $m_{r_0j_0l}\neq 0$ for some $l$, the polynomial
$R_{r_0j_0-1}$ is nonconstant. We have thus established uniform
distribution for a set of $x$ of full measure, completing the
proof.
\end{proof}

\noindent{\bf Acknowledgment}:  The authors thank the referee for
his help in organizing and simplifying the presentation, and in particular
for the simple proof of Proposition~\ref{C:1}.

\end{document}